\documentclass{amsart}

\usepackage{amsmath} 
\usepackage{amssymb}
\usepackage{mathrsfs}

\newtheorem{theorem}{Theorem}[section] 
\newtheorem{claim}[theorem]{Claim}

\newtheorem{conclusion}[theorem]{Conclusion}
\newtheorem{observation}[theorem]{Observation}

\theoremstyle{definition}
\newtheorem{definition}[theorem]{Definition}

\theoremstyle{remark}

\newcommand{\pp}{{\rm pp}}

\newcommand{\ch}{{\rm ch}}

\newcommand{\Ord}{{\rm Ord}}

\newcommand{\INC}{{\rm INC}}

\newcommand{\chr}{{\rm chr}}
\newcommand{\otp}{{\rm otp}}

\newcommand{\edge}{{\rm edge}}

\newcommand{\nacc}{{\rm nacc}}

\newcommand{\Dom}{{\rm Dom}}
\newcommand{\Rang}{{\rm Rang}}

\newcommand{\rest}{{\restriction}}

\newcommand{\wilog}{{\rm without loss of generality}}
\newcommand{\Wilog}{{\rm Without loss of generality}}

\newcommand{\then}{{\underline{then}}}
\newcommand{\when}{{\underline{when}}}
\newcommand{\Then}{{\underline{Then}}}

\newcommand{\Iff}{{\underline{iff}}}
\newcommand{\mn}{{\medskip\noindent}}
\newcommand{\sn}{{\smallskip\noindent}}

\newcommand{\cA}{{\mathscr A}}

\newcommand{\cB}{{\mathscr B}}

\newcommand{\cf}{{\rm cf}}

\newcount\skewfactor
\def\mathunderaccent#1#2 {\let\theaccent#1\skewfactor#2
\mathpalette\putaccentunder}
\def\putaccentunder#1#2{\oalign{$#1#2$\crcr\hidewidth
\vbox to.2ex{\hbox{$#1\skew\skewfactor\theaccent{}$}\vss}\hidewidth}}

\newenvironment{PROOF}[2][\proofname.]
   {\begin{proof}[#1]\renewcommand{\qedsymbol}{\textsquare$_{\rm #2}$}}
   {\end{proof}}

\begin{document}

\title {On incompactness for chromatic number of graphs}
\author {Saharon Shelah}
\address{Einstein Institute of Mathematics\\
Edmond J. Safra Campus, Givat Ram\\
The Hebrew University of Jerusalem\\
Jerusalem, 91904, Israel\\
 and \\
 Department of Mathematics\\
 Hill Center - Busch Campus \\ 
 Rutgers, The State University of New Jersey \\
 110 Frelinghuysen Road \\
 Piscataway, NJ 08854-8019 USA}
\email{shelah@math.huji.ac.il}
\urladdr{http://shelah.logic.at}
\thanks{The author thanks Alice Leonhardt for the beautiful typing.\\
The author would like to thank the Israel Science Foundation for partial
support of this research (Grant no. 1053/11). Publication 1006.}



\date{June 28, 2012}

\begin{abstract}
We deal with incompactness.  Assume the existence of non-reflecting
 stationary set of cofinality $\kappa$.  We prove that one can
 define a graph $G$ whose chromatic number is $>\kappa$, while the
 chromatic number of every subgraph $G' \subseteq G,|G'| < |G|$ is
 $\le \kappa$.  The main case is $\kappa = \aleph_0$.
\end{abstract}

\maketitle
\numberwithin{equation}{section}
\setcounter{section}{-1}
\newpage

\centerline {Anotated Content}
\bigskip

\noindent
\S0 \quad Introduction, pg.\pageref{Introduction}

\S(0A) \quad The questions and results, pg.\pageref{Thequestions}

\S(0AB) \quad Preliminaries, pg.\pageref{Preliminaries}
\bigskip

\noindent
\S1 \quad From a non-reflecting stationary set, pg.\pageref{From}
\mn
\begin{enumerate}
\item[${{}}$]  [We show that ``$S \subseteq S^\lambda_\kappa$ is
  stationary not reflecting" implies imcompactness for length
  $\lambda$ for ``chromatic number $=\kappa$".]
\end{enumerate}
\bigskip

\noindent
\S2 \quad From almost free, pg.\pageref{Fromalmost}
\mn
\begin{enumerate}
\item[${{}}$]  [Here we weaken the assumption in \S1 to ``$\cA
  \subseteq {}^{\kappa}\Ord$ is almost free".]
\end{enumerate}
\newpage

\section {Introduction} \label{Introduction}

\subsection {The questions and results} \label{Thequestions}
\bigskip

During the Hajnal conference (June 2011) Magidor asked me on
incompactness of ``having chromatic number $\aleph_0$"; that is,
there is a graph $G$ with $\lambda$ nodes, chromatic number $>
\aleph_0$ but every subgraph with $< \lambda$ nodes has chromatic
number $\aleph_0$ when:
\mn
\begin{enumerate}
\item[$(*)_1$]  $\lambda$ is regular $> \aleph_1$ with a
  non-reflecting stationary $S \subseteq S^\lambda_{\aleph_0}$,
  possibly though better not, assuming some version of GCH.
\end{enumerate}
\mn
Subsequently also when:
\mn
\begin{enumerate}
\item[$(*)_2$]  $\lambda = \aleph_{\omega +1}$.
\end{enumerate}
\mn
Such problems were first asked by Erd\"os-Hajnal, see \cite{EH1}; we
continue \cite{Sh:347}.

First answer was using BB, see \cite[3.24]{Sh:309} so assuming
\mn
\begin{enumerate}
\item[$\boxplus$]  $(a) \quad \lambda = \mu^+$
\sn
\item[${{}}$]  $(b) \quad \mu^{\aleph_0} = \mu$
\sn
\item[${{}}$]  $(c) \quad S \subseteq \{\delta < \lambda:\cf(\delta) =
  \aleph_0\}$ is stationary not reflecting
\end{enumerate}
\mn
or just
\mn
\begin{enumerate}
\item[$\boxplus'$]  $(a) \quad \lambda = \cf(\lambda)$
\sn
\item[${{}}$]  $(b) \quad \alpha < \lambda \Rightarrow
  |\alpha|^{\aleph_0} < \lambda$
\sn
\item[${{}}$]  $(c) \quad$ as above.
\end{enumerate}
\mn
However, eventually we get more: if $\lambda = \lambda^{\aleph_0} =
\cf(\lambda)$ and $S \subseteq S^\lambda_{\aleph_0}$ is stationary
non-reflective then we have $\lambda$-incompactness for
$\aleph_0$-chromatic.  In fact,  we replace $\aleph_0$ by
$\kappa = \cf(\kappa) < \lambda$ using a suitable hypothesis.

Moreover, if $\lambda^\kappa > \lambda$ we still get
$(\lambda^\kappa,\lambda)$-incompactness for $\kappa$-chromatic
number.
In \S2 we use quite free family of countable sequences.

In subsequent work we shall solve also the parallel of the second
question of Magidor, i.e. 
\mn
\begin{enumerate}
\item[$(*)_2$]  for regular $\kappa \ge \aleph_0$ and
 $\varepsilon < \kappa$ there is a graph $G$ of chromatic number $>
 \kappa$ but every sub-graph with $< \aleph_{\kappa \cdot \varepsilon
 +1}$ nodes has chromatic number $\le \kappa$.
\end{enumerate}
\mn
We thank Menachem Magidor for asking, Peter Komjath for stimulating
discussion and Paul Larson, Shimoni Garti and the referee for some comments.
\bigskip

\subsection{Preliminaries} \label{Preliminaries} \
\bigskip

\begin{definition}
\label{y3}
For a graph $G$, let $\ch(G)$, the chromatic number of $G$ be the
minimal cardinal $\chi$ such that there is colouring $\bold c$ of $G$
with $\chi$ colours, that is $\bold c$ is a function from the set of
nodes of $G$ into $\chi$ or just a set of of cardinality $\le \chi$ such
that $\bold c(x) = \bold c(y) \Rightarrow \{x,y\} \notin \edge(G)$.
\end{definition}

\begin{definition}
\label{y6}
1) We say ``we have $\lambda$-incompactness for the $(< \chi)$-chromatic
   number" or $\INC_{\chr}(\lambda,< \chi)$ \when \,: 
there is a graph $G$ with $\lambda$ nodes, chromatic number $\ge \chi$
   but every subgraph with $< \lambda$ nodes has
   chromatic number $< \chi$.

\noindent
2) If $\chi = \mu^+$ we may replace $``< \chi"$ by $\mu$; 
similarly in \ref{y8}.
\end{definition}

\noindent
We also consider
\begin{definition}
\label{y8}
1) We say ``we have $(\mu,\lambda)$-incompactness for $(< \chi)$-chromatic
 number" or $\INC_{\chr}(\mu,\lambda,< \chi)$
 \when \, there is an increasing continuous sequence $\langle
 G_i:i \le \lambda\rangle$ of graphs each with $\le \mu$ nodes, $G_i$ an
 induced subgraph of $G_\lambda$ with $\ch(G_\lambda) \ge \chi$ but $i
 < \lambda \Rightarrow \ch(G_i) < \chi$.

\noindent
2) Replacing (in part (1)) $\chi$ by $\bar\chi = (< \chi_0,\chi_1)$
means $\ch(G_\lambda)) \ge \chi_1$ and $i < \lambda \rightarrow \ch(G_i) <
\chi_0$; similarly in \ref{y6} and parts 3),4) below.

\noindent
3) We say we have incompactness for length $\lambda$ for $(< \chi)$-chromatic
(or $\bar\chi$-chromatic) number \when \, we fail to have
$(\mu,\lambda)$-compactness for $(< \chi)$-chromatic (or
   $\bar\chi$-chromatic) number for some $\mu$.

\noindent
4) We say we have $[\mu,\lambda]$-incompactness for $(<
   \chi)$-chromatic number or $\INC_{\chr}[\mu,\lambda,< \chi]$ 
\when \, there is a graph $G$ with $\mu$ nodes, 
$\ch(G) \ge \chi$ but $G^1 \subseteq G \wedge |G^1| < \lambda \Rightarrow
   \ch(G^1) < \chi$.

\noindent
5) Let $\INC^+_{\chr}(\mu,\lambda,< \chi)$ be as in part (1) but we
   add that even the $c \ell(G_i)$, the colouring number of $G_i$ is $<
   \chi$ for $i < \lambda$, see below.

\noindent
6) Let $\INC^+_{\chr}[\mu,\lambda,< \chi]$ be as in part (4) but we
add $G^1 \subseteq G \wedge |G^1| < \lambda \Rightarrow c \ell(G^1) <
\chi$.

\noindent
7) If $\chi = \kappa^+$ we may write $\kappa$ instead of ``$< \chi$".
\end{definition}

\begin{definition}
\label{y11}
1) For regular $\lambda > \kappa$ let $S^\lambda_\kappa = \{\delta <
   \lambda:\cf(\delta) = \kappa\}$.

\noindent
2) We say $C$ is a $(\ge \theta)$-closed subset of a set $B$ of
   ordinals when: if $\delta = \sup(\delta \cap B) \in B,\cf(\delta)
   \ge \theta$ and $\delta = \sup(C \cap \delta)$ then $\delta \in C$.
\end{definition}

\begin{definition}
\label{y13}
For a graph $G$, the colouring number $c \ell(G)$ is the minimal
$\kappa$ such that there is a list $\langle a_\alpha:\alpha <
\alpha(*)\rangle$ of the nodes of $G$ such that $\alpha < \alpha(*)
\Rightarrow \kappa > |\{\beta < \alpha:\{a_\beta,a_\alpha\} \in \edge(G)\}$.
\end{definition}
\newpage

\section {From non-reflecting stationary in cofinality $\aleph_0$} \label{From}

\begin{claim}  
\label{a3}
There is a graph $G$ with $\lambda$ nodes and chromatic number $>
\kappa$ but every subgraph with $< \lambda$ nodes have
chromatic number $\le \kappa$ \when \,:
\mn
\begin{enumerate}
\item[$\boxplus$]  $(a) \quad \lambda,\kappa$ are regular cardinals
\sn
\item[${{}}$]  $(b) \quad \kappa < \lambda = \lambda^\kappa$
\sn
\item[${{}}$]  $(c) \quad S \subseteq S^\lambda_\kappa$ 
is stationary, not reflecting.
\end{enumerate}
\end{claim}

\begin{PROOF}{\ref{a3}}

\noindent
\underline{Stage A}:  Let $\bar X = \langle X_i:i < \lambda\rangle$ be a
partition of $\lambda$ to sets such that $|X_i| = \lambda$ or just $|X_i|
= |i+2|^\kappa$ and $\min(X_i) \ge i$ and let $X_{<i} = \cup\{X_j:j
< i\}$ and $X_{\le i} = X_{<(i+1)}$.
For $\alpha < \lambda$ let $\bold i(\alpha)$ be the unique ordinal
$i < \lambda$ such that 
$\alpha \in X_i$.  We choose the set of points = nodes of
$G$ as $Y = \{(\alpha,\beta):\alpha < \beta < \lambda,\bold i(\beta) \in
S$ and $\alpha < \bold i(\beta)\}$ and let $Y_{<i} = \{(\alpha,\beta) \in
Y:\bold i(\beta) < i\}$.
\medskip

\noindent
\underline{Stage B}:  Note that if $\lambda = \kappa^+$, the complete
graph with $\lambda$ nodes is an example (no use of the further
information in $\boxplus$).  So \wilog \, $\lambda > \kappa^+$.

Now choose a sequence satisfying the following properties, 
exists by \cite[Ch.III]{Sh:g}:
\mn
\begin{enumerate}
\item[$\boxplus$]  $(a) \quad \bar C = \langle C_\delta:\delta \in
  S\rangle$
\sn
\item[${{}}$]  $(b) \quad C_\delta \subseteq \delta = \sup(C_\delta)$
\sn
\item[${{}}$]  $(c) \quad \otp(C_\delta) = \kappa$ such that $(\forall
  \beta \in C_\delta)(\beta +1,\beta +2 \notin C_\delta)$
\sn
\item[${{}}$]  $(d) \quad \bar C$ guesses\footnote{the guessing clubs
  are used only in Stage D.}clubs.
\end{enumerate}
\mn
Let $\langle \alpha^*_{\delta,\varepsilon}:\varepsilon <
\kappa\rangle$ list $C_\delta$ in increasing order.

For $\delta \in S$ let $\Gamma_\delta$ be the set of sequence
$\bar\beta$ such that:
\mn
\begin{enumerate}
\item[$\boxplus_{\bar\beta}$]  $(a) \quad \bar \beta$ has the form
$\langle \beta_\varepsilon:\varepsilon < \kappa \rangle$
\sn
\item[${{}}$]  $(b) \quad \bar\beta$ is increasing with limit $\delta$
\sn
\item[${{}}$]  $(c) \quad \alpha^*_{\delta,\varepsilon} < 
\beta_{2 \varepsilon +i} < \alpha^*_{\delta,\varepsilon +1}$ for $i <
2,\varepsilon < \kappa$
\sn
\item[${{}}$]  $(d) \quad \beta_{2 \varepsilon +i} \in 
X_{< \alpha^*_{\delta,\varepsilon +1}} \backslash X_{\le
  \alpha^*_{\delta,\varepsilon}}$ for $i < 2,\varepsilon < \kappa$
\sn
\item[${{}}$]  $(e) \quad (\beta_{2 \varepsilon},\beta_{2\varepsilon
  +1}) \in Y$ hence $\in Y_{< \alpha^*_{\delta,\varepsilon +1}} \subseteq
Y_{< \delta}$ for each $\varepsilon < \kappa$
\end{enumerate}
\mn
(can ask less).

So $|\Gamma_\delta| \le |\delta|^\kappa \le |X_\delta| \le \lambda$ 
hence we can choose a sequence
$\langle \bar\beta_\gamma:\gamma \in X'_\delta \subseteq
X_\delta\rangle$ listing $\Gamma_\delta$.

Now we define the set of edges of $G$: $\edge(G) =
\{\{(\alpha_1,\alpha_2),(\min(C_\delta),\gamma)\}:\delta \in S,\gamma
\in X'_\delta$ hence the sequence $\bar\beta_\gamma = \langle
\beta_{\gamma,\varepsilon}:\varepsilon < \kappa\rangle$ is well
defined and we demand $(\alpha_1,\alpha_2) \in \{(\beta_{\gamma,2
  \varepsilon},\beta_{\gamma,2 \varepsilon +1}):\varepsilon < \kappa\}\}$.
\medskip

\noindent
\underline{Stage C}:  Every subgraph of $G$ of cardinality $<
 \lambda$ has chromatic number $\le \kappa$.

For this we shall prove that:
\mn
\begin{enumerate}
\item[$\oplus_1$]  $\ch(G \rest Y_{<i}) \le \kappa$ for every $i <
  \lambda$.
\end{enumerate}
\mn
This suffice as $\lambda$ is regular, hence every subgraph with $<
\lambda$ nodes is included in $Y_{<i}$ for some $i < \lambda$.

For this we shall prove more by induction on $j < \lambda$:
\mn
\begin{enumerate}
\item[$\oplus_{2,j}$]  if $i < j,i \notin S,\bold c_1$ a colouring of
  $G \rest Y_{<i},\Rang(\bold c_1) \subseteq \kappa$ and $u \in
[\kappa]^\kappa$ \then \, there is a colouring $\bold c_2$ of 
$G \rest Y_{<j}$ extending $\bold c_1$ such that $\Rang(\bold c_2 \rest
  (Y_{<j} \backslash Y_{<i})) \subseteq u$.
\end{enumerate}
\medskip

\noindent
\underline{Case 1}:  $j=0$

Trivial.
\medskip

\noindent
\underline{Case 2}:  $j$ successor, $j-1 \notin S$

By the induction hypothesis \wilog \, $j=i+1$, but then every node
from $Y_j \backslash Y_i$ is an isolated node in $G \rest Y_{<j}$,
because if $\{(\alpha,\beta),(\alpha',\beta')\}$ is an edge of $G
\rest Y_j$ then $\bold i(\beta),\bold i(\beta') \in S$ hence
necessarily $\bold i(\beta) \ne j-1=i,\bold i(\beta') \ne j-1=i$ hence
both $(\alpha,\beta),(\alpha,\beta')$ are from $Y_i$.
\medskip

\noindent
\underline{Case 3}:  $j$ successor, $j-1 \in S$

Let $j-1$ be called $\delta$ so $\delta \in S$.  But $i \notin S$ 
by the assumption in $\oplus_{2,j}$ hence $i < \delta$.  
Let $\varepsilon(*) < \kappa$ be such that
$\alpha^*_{\delta,\varepsilon(*)} > i$.

Let $\langle u_\varepsilon:\varepsilon \le \kappa\rangle$ be a
sequence of subsets of $u$, a partition of $u$ to sets each of 
cardinality $\kappa$; actually the only disjointness used is
that $u_\kappa \cap (\bigcup\limits_{\varepsilon < \kappa}
u_\varepsilon) = \emptyset$.

We let $i_0 = i,i_{1 + \varepsilon} =
\cup\{\alpha^*_{\delta,\varepsilon(*)+1+\zeta} +1:\zeta < 1 +
\varepsilon\},i_\kappa = \delta,i_{\kappa +1} = \delta +1 =j$.

Note that:
\mn
\begin{enumerate}
\item[$\bullet$]  $\varepsilon < \kappa \Rightarrow i_\varepsilon
\notin S_j$.
\end{enumerate}
\mn
[Why?  For $\varepsilon=0$ by the assumption on $i$, for $\varepsilon$
successor $i_\varepsilon$ is a successor ordinal and for $i$ limit clearly
$\cf(i_\varepsilon) = \cf(\varepsilon) < \kappa$ and $S \subseteq
S^\lambda_\kappa$.]

We now choose $\bold c_{2,\zeta}$ by induction on $\zeta \le \kappa +1$ such
that:
\mn
\begin{enumerate}
\item[$\bullet$]  $\bold c_{2,0} = \bold c_1$
\sn
\item[$\bullet$]  $\bold c_{2,\zeta}$ is a colouring of $G \rest Y_{<i_\zeta}$
\sn
\item[$\bullet$]  $\bold c_{2,\zeta}$ is increasing with $\zeta$
\sn
\item[$\bullet$]  $\Rang(\bold c_{2,\zeta} \rest (Y_{< i_{\xi +1}}
  \backslash Y_{< i_\xi})) \subseteq u_\xi$ for every $\xi < \zeta$.
\end{enumerate}
\mn
For $\zeta = 0,\bold c_{2,0}$ is $\bold c_1$ so is given.

For $\zeta = \varepsilon +1 < \kappa$: use the induction hypothesis,
possible as necessarily $i_\varepsilon \notin S$.

For $\zeta \le \kappa$ limit: take union.

For $\zeta = \kappa+1$, note that each node $b$ of $Y_{< i_\zeta}
\backslash Y_{< i_\kappa}$ is not connected to any other such node and
if the node $b$ is connected to a node from $Y_{< i_\kappa}$ then the
node $b$ necessarily has the form
$(\min(C_\delta),\gamma),\gamma \in X'_\delta$, hence $\bar\beta_\gamma$
is well defined, so the node $b = (\min(C_\delta),\gamma)$
is connected in $G$, more exactly in $G \rest Y_{\le \delta}$ 
exactly to the $\kappa$ nodes
$\{(\beta_{\gamma,2 \varepsilon},\beta_{\gamma,2 \varepsilon
  +1}):\varepsilon < \kappa\}$, but for every $\varepsilon < \kappa$
large enough, $\bold c_{2,\kappa}((\beta_{\gamma,2
\varepsilon},\beta_{\gamma,2 \varepsilon +1})) \in 
u_\varepsilon$ hence $\notin u_\kappa$ and
$|u_\kappa| = \kappa$ so we can choose a colour.
\medskip

\noindent
\underline{Case 4}:  $j$ limit

By the assumption of the claim there is a club $e$ 
of $j$ disjoint to $S$ and \wilog \, $\min(e) =
i$.  Now choose $\bold c_{2,\xi}$ a colouring of $Y_{< \xi}$ 
by induction on $\xi \in e \cup \{j\}$, increasing with $\xi$ such that
$\Rang(\bold c_{2,\xi} \rest (Y_{< \varepsilon} \backslash Y_{<i}))
\subseteq u$ and $\bold c_{2,0} = \bold c_1$
\mn
\begin{enumerate}
\item[$\bullet$]   For $\xi = \min(e) = i$ the colouring $\bold
 c_{2,\xi} = \bold c_{2,i} = \bold c_1$ is given,
\sn
\item[$\bullet$]   for $\xi$ successor in $e$, i.e. $\in
\nacc(e) \backslash \{i\}$, use the induction hypothesis with
$\xi,\max(e \cap \xi)$ here playing the role of $j,i$ there recalling
$\max(e \cap \xi) \in e,e \cap S = \emptyset$
\sn
\item[$\bullet$]  for $\xi = \sup(e \cap \xi)$ take union.
\end{enumerate}
\mn
Lastly, for $\xi=j$ we are done.
\medskip

\noindent
\underline{Stage D}:  $\ch(G) > \kappa$.

Why?  Toward a contradiction, assume $\bold c$ is 
a colouring of $G$ with set of colours
$\subseteq \kappa$.  For each $\gamma < \lambda$ let 
$u_\gamma = \{\bold c((\alpha,\beta)):\gamma < \alpha <
\beta < \lambda$ and $(\alpha,\beta) \in Y\}$.  So $\langle u_\gamma:\gamma 
< \lambda\rangle$ is $\subseteq$-decreasing sequence of subsets of $\kappa$ and
$\kappa < \lambda = \cf(\lambda)$, 
hence for some $\gamma(*) < \lambda$ and $u_*
\subseteq \kappa$ we have $\gamma \in (\gamma(*),\lambda) 
\Rightarrow u_\gamma = u_*$.

Hence $E = \{\delta < \lambda:\delta$ is a limit ordinal $> \gamma(*)$
and $(\forall \alpha < \delta)((\bold i(\alpha) < \delta)$ 
and for every $\gamma < \delta$ and $i \in u_*$ there are $\alpha <
\beta$ from $(\gamma,\delta)$ such that $(\alpha,\beta) \in Y$ and
$\bold c((\alpha,\beta))=i\}$ is a club of $\lambda$.

Now recall that $\bar C$ guesses clubs hence for some $\delta \in S$ we have
$C_\delta \subseteq E$, so for every $\varepsilon < \kappa$ we can
choose $\beta_{2 \varepsilon} < \beta_{2 \varepsilon +1}$ from
$(\alpha^*_{\delta,\varepsilon},\alpha^*_{\delta,\varepsilon +1})$
such that $(\beta_{2 \varepsilon},\beta_{2 \varepsilon +1}) \in Y$ and
$\varepsilon \in u_* \Rightarrow \bold c((\beta_{2
  \varepsilon},\beta_{2 \varepsilon +1})) = \varepsilon$.  So $\langle
\beta_\varepsilon:\varepsilon < \kappa\rangle$ is well defined, increasing and
belongs to $\Gamma_\delta$, hence $\bar\beta_\gamma = \langle
\beta_\varepsilon:\varepsilon < \kappa\rangle$ for some $\gamma \in
X_\delta$, hence $(\alpha^*_{\delta,0},\gamma)$ belongs to $Y$ and is
connected in the graph to $(\beta_{2 \varepsilon},\beta_{2 \varepsilon
  +1})$ for $\varepsilon < \kappa$.  Now if $\varepsilon\in u_*$ then
$\bold c((\beta_{2 \varepsilon},\beta_{2 \varepsilon +1})) =
\varepsilon$ hence $\bold c((\alpha^*_{\delta,0},\gamma)) \ne
\varepsilon$ for every $\varepsilon \in u_*$, 
so $\bold c((\alpha^*_{\delta,0},\gamma)) \in \kappa
\backslash u_*$.  But $u_* = u_{\alpha^*_{\delta,0}}$ and $\bold
c((\alpha^*_{\delta,0},\gamma)) \in \kappa \backslash u_*$, so we get
contradiction to the definition of $u_{\alpha^*_{\delta,0}}$.
\end{PROOF}

\noindent
Similarly
\begin{claim}
\label{a6}
There is an increasing continuous sequence $\langle G_i:i \le
\lambda \rangle$ of graphs each of cardinality $\lambda^\kappa$ such
that $\ch(G_\lambda) > \kappa$ and $i < \lambda$ implies $\ch(G_i) \le
\kappa$ and even $c \ell(G_i) \le \kappa$ \when \,:
\mn
\begin{enumerate}
\item[$\boxplus$]  $(a) \quad \lambda = \cf(\lambda)$
\sn
\item[${{}}$]  $(b) \quad S \subseteq \{\delta < \lambda:\cf(\delta) =
  \kappa\}$ is stationary not reflecting.
\end{enumerate}
\end{claim}

\begin{PROOF}{\ref{a6}}
Like \ref{a3} but the $X_i$ are not necessarily $\subseteq \lambda$ or
 use \ref{c3}.
\end{PROOF}
\newpage

\section {From almost free} \label{Fromalmost}

\begin{definition}
\label{c1}
Suppose $\eta_\beta \in {}^\kappa\Ord$ for every $\beta < \alpha(*)$
and $u \subseteq \alpha(*)$, and $\alpha < \beta < \alpha(*)
\Rightarrow \eta_\alpha \ne \eta_\beta$.

\noindent
1) We say $\{\eta_\alpha:\alpha \in u\}$ is free \when \, there exists a
function $h:u \rightarrow \kappa$ such that $\langle
   \{\eta_\alpha(\varepsilon):\varepsilon \in
   [h(\alpha),\kappa)\}:\alpha \in u\rangle$ is a sequence of pairwise
   disjoint sets.

\noindent
2) We say $\{\eta_\alpha:\alpha \in u\}$ is weakly free \when \, there
exists a sequence $\langle u_{\varepsilon,\zeta}:\varepsilon,\zeta <
   \kappa\rangle$ of subsets of $u$ with union $u$, such that the
   function $\eta_\zeta \mapsto \eta_\zeta(\varepsilon)$ is a
   one-to-one function on $u_{\varepsilon,\zeta}$, 
for each $\varepsilon,\zeta < \kappa$.
\end{definition}

\begin{claim}
\label{c3}
1) We have $\INC_{\chr}(\mu,\lambda,\kappa)$ and even
$\INC^+_{\chr}(\mu,\lambda,\kappa)$, see Definition \ref{y8}(1),(5) \when \,:
\mn
\begin{enumerate}
\item[$\boxplus$]  $(a) \quad \alpha(*) \in [\mu,\mu^+)$ and $\lambda$
  is regular $\le \mu$ and $\mu = \mu^\kappa$
\sn
\item[${{}}$]  $(b) \quad \bar \eta = \langle \eta_\alpha:\alpha <
\alpha(*)\rangle$ 
\sn
\item[${{}}$]  $(c) \quad \eta_\alpha \in {}^\kappa \mu$
\sn
\item[${{}}$]  $(d) \quad \langle u_i:i \le \lambda\rangle$
is a $\subseteq$-increasing continuous sequence of subsets of
$\alpha(*)$

\hskip25pt  with $u_\lambda = \alpha(*)$
\sn
\item[${{}}$]  $(e) \quad \bar \eta \rest u_\alpha$ is free \Iff \,
$\alpha < \lambda$ \Iff \, $\bar\eta \rest u_\alpha$ is weakly free.
\end{enumerate}
\mn
2) We have $\INC_{\chr}[\mu,\lambda,\kappa]$ and even
$\INC^+_{\chr}[\mu,\lambda,\kappa]$ , see Definition \ref{y8}(4) \when \,:
\mn
\begin{enumerate}
\item[$\boxplus_2$]  $(a),(b),(c) \quad$ as in $\boxplus$ from \ref{c3}
\sn
\item[${{}}$]  $(d) \quad \bar\eta$ is not free
\sn
\item[${{}}$]  $(e) \quad \bar\eta \rest u$ is free when 
$u \in [\alpha(*)]^{< \lambda}$. 
\end{enumerate}
\end{claim}

\begin{PROOF}{\ref{c3}}
We concentrate on proving part (1); the proof of part (2) is similar.
For $\cA \subseteq {}^\kappa\Ord$, we define $\tau_{\cA}$ as the
vocabulary $\{P_\eta:\eta \in \cA\} \cup \{F_\varepsilon:\varepsilon <
\kappa\}$ where $P_\eta$ is a unary predicate, $F_\varepsilon$ a unary
function (will be interpreted as possibly partial).

\Wilog \, for each $i < \lambda,u_i$ is an initial
segment of $\alpha(*)$ and let 
$\cA = \{\eta_\alpha:\alpha < \alpha(*)\}$
and let $<_{\cA}$ be the well ordering 
$\{(\eta_\alpha,\eta_\beta):\alpha < \beta < \alpha(*)\}$ of $\cA$.

We further let $K_{\cA}$ be the class of structures $M$ such that
(pedantically, $K_{\cA}$ depend also on the sequence $\langle
\eta_\alpha:\alpha < \alpha(*)\rangle$:
\mn
\begin{enumerate}
\item[$\boxplus_1$]  $(a) \quad M = (|M|,F^M_\varepsilon,
P^M_\eta)_{\varepsilon < \kappa,\eta \in \cA}$
\sn
\item[${{}}$]  $(b) \quad \langle P^M_\eta:\eta \in \cA\rangle$ is a
partition of $|M|$, so for $a \in M$ let $\eta_a$ 

\hskip25pt $= \eta^M_a$ be the unique $\eta \in \cA$ such that $a \in P^M_\eta$
\sn
\item[${{}}$]  $(c) \quad$ if $a_\ell \in P^M_{\eta_\ell}$ for
  $\ell=1,2$ and $F^M_\varepsilon(a_2) = a_1$ then

\hskip25pt $\eta_1(\varepsilon) = \eta_2(\varepsilon)$ and
$\eta_1 <_{\cA} \eta_2$.
\end{enumerate}
\mn
Let $K^*_{\cA}$ be the class of $M$ such that
\mn
\begin{enumerate}
\item[$\boxplus_2$]  $(a) \quad M \in K_{\cA}$
\sn
\item[${{}}$]  $(b) \quad \|M\| = \mu$
\sn
\item[${{}}$]  $(c) \quad$ if $\eta \in \cA,u \subseteq \kappa$ and
$\eta_\varepsilon <_{\cA} \eta,\eta_\varepsilon(\varepsilon) 
= \eta(\varepsilon)$ and $a_\varepsilon \in P^M_{\eta_\varepsilon}$ 

\hskip25pt  for $\varepsilon \in u$ \then \, for some $a \in P^M_\eta$ we have 
$\varepsilon \in u \Rightarrow F^M_\varepsilon(a) = a_\varepsilon$ 

\hskip25pt and $\varepsilon \in \kappa
  \backslash  u \Rightarrow F^M_\varepsilon(a)$ not defined.
\end{enumerate}
\mn
Clearly
\mn
\begin{enumerate}
\item[$\boxplus_3$]   there is $M \in K^*_{\cA}$.
\end{enumerate}
\mn
[Why?  As $\mu = \mu^\kappa$ and $|\cA| = \mu$.]
\mn
\begin{enumerate}
\item[$\boxplus_4$]  for $M \in K_{\cA}$ let $G_M$ be the graph with:
\sn
\begin{enumerate}
\item[$\bullet$]  set of nodes $|M|$
\sn
\item[$\bullet$]  set of edges $\{\{a,F^M_\varepsilon(a)\}:a \in
  |M|,\varepsilon < \kappa$ when $F^M_\varepsilon(a)$ is defined$\}$.
\end{enumerate}
\end{enumerate}
\mn
Now
\mn
\begin{enumerate}
\item[$\boxplus_5$]  if $u \subseteq \alpha(*),\cA_u =
\{\eta_\alpha:\alpha \in u\} \subseteq \cA$ and $\bar\eta \rest u$
is free, and $M \in K_{\cA}$ \then \, $G_{M,\cA_u} := G_M \rest 
(\cup\{P^M_\eta:\eta \in \cA_u\})$
has chromatic number $\le \kappa$; moreover has colouring number $\le \kappa$.
\end{enumerate}
\mn
[Why?  Let $h:u \rightarrow \kappa$ witness that $\bar\eta \rest u$ is
free and for $\varepsilon < \kappa$ let $\cB_\varepsilon :=
\{\eta_\alpha:\alpha \in u$ and $h(\alpha) = \varepsilon\}$, so $\cB =
\cup\{\cB_\varepsilon:\varepsilon < \kappa\}$, hence it is enough to
prove for each $\varepsilon < \kappa$ that $G_{\mu,\cB_\varepsilon}$
has chromatic number $\le \kappa$.  To prove this, by induction on
$\alpha \le \alpha(*)$ we choose $\bold c^\varepsilon_\alpha$ such that:
\mn
\begin{enumerate}
\item[$\boxplus_{5.1}$]  $(a) \quad \bold c^\varepsilon_\alpha$ is a
  function 
\sn
\item[${{}}$]  $(b) \quad \langle \bold c_\beta:\beta \le
  \alpha\rangle$ is increasing continuous
\sn
\item[${{}}$]  $(c) \quad \Dom(\bold c^\varepsilon_\alpha) =
  B^\varepsilon_\alpha := \cup\{P^M_{\eta_\beta}:\beta < \alpha$ and
  $\eta_\beta \in \cB_\varepsilon\}$
\sn
\item[${{}}$]  $(d) \quad \Rang(\bold c^\varepsilon_\alpha) \subseteq
  \kappa$
\sn
\item[${{}}$]  $(e) \quad$ if $a,b, \in \Dom(\bold c_\alpha)$ and
  $\{a,b\} \in \edge(G_M)$ then $\bold c_\alpha(a) \ne \bold
  c_\alpha(b)$.
\end{enumerate}
\mn
Clearly this suffices.  Why is this possible?

If $\alpha = 0$ let $\bold c^\varepsilon_\alpha$ be empty, if $\alpha$
is a limit ordinal let $\bold c^\varepsilon_\alpha = \cup\{\bold
c^\varepsilon_\beta:\beta < \alpha\}$ and if $\alpha = \beta +1 \wedge
\alpha(\beta) \ne G$ let $\bold c_\alpha = \bold c_\beta$.

Lastly, if $\alpha = \beta +1 \wedge h(\beta) = \varepsilon$ we define
$\bold c^\varepsilon_\alpha$ as follows for $a \in \Dom(\bold
c^\varepsilon_\alpha),\bold c^\varepsilon_\alpha(a)$ is:
\bigskip

\noindent
\underline{Case 1}:  $a \in B^\varepsilon_\beta$.

Then $\bold c^\varepsilon_\alpha(a) = \bold c^\varepsilon_\beta(a)$.
\bigskip

\noindent
\underline{Case 2}:  $a \in B^\varepsilon_\alpha \backslash
B^\varepsilon_\beta$.

Then $\bold c^\varepsilon_\alpha(a) = \min(\kappa
\backslash \{\bold c^\varepsilon_\beta(F^M_\zeta(a)):\zeta < \varepsilon$
and $F^M_\zeta(a) \in \Dom(\bold c^\varepsilon_\beta)\})$.

This is well defined as:
\mn
\begin{enumerate}
\item[$\boxplus_{5.2}$]  $(a) \quad B^\varepsilon_\alpha =
  B^\varepsilon_\beta \cup P^M_{\eta_\beta}$
\sn
\item[${{}}$]  $(b) \quad$ if $a \in B^\varepsilon_\beta$ then $\bold
  c^\varepsilon_\beta(a)$ is well defined (so case 1 is O.K.)
\sn
\item[${{}}$]  $(c) \quad$ if $\{a,b\} \in \edge(G_M),a \in
P^M_{\eta_\beta}$ and $b \in B^\varepsilon_\alpha$ \then \, $b \in
  B^\varepsilon_\beta$ and

\hskip25pt  $b \in \{F^M_\zeta(a):\zeta < \varepsilon\}$
\sn
\item[${{}}$]  $(d) \quad \bold c^\varepsilon_\alpha(a)$ is well
  defined in Case 2, too
\sn
\item[${{}}$]  $(e) \quad  \bold c^\varepsilon_\alpha$ is a function
  from $B^\varepsilon_\alpha$ to $\kappa$
\sn
\item[${{}}$]  $(f) \quad \bold c^\varepsilon_\alpha$ is a colouring.
\end{enumerate}
\mn
[Why?  Clause (a) by $\boxplus_{5.1}(c)$, clause (b) by the induction
  hypothesis and clause (c) by $\boxplus_1(c) + \boxplus_4$.  Next,
  clause (d) holds as $\{\bold c^\varepsilon_\beta(F^M_\zeta(a)):\zeta <
  \varepsilon$ and $F^M_\zeta(a) \in B^\varepsilon_\beta = \Dom(\bold
  c^\varepsilon_\beta )\}$ is a set of cardinality $\le |\varepsilon|
  < \kappa$.  Clause (e) holds by the choices of the $\bold
  c^\varepsilon_\alpha(a)$'s.  Lastly, to check that clause (f) holds
assume $(a,b)$ is an edge of $G_M \rest B^\varepsilon_\alpha$, for
some $\zeta < \kappa$ we have $b = F^M_\zeta(a)$, hence $\eta^M_a <_{\cA}
\eta^M_b$.  If $a,b \in B^\varepsilon_\beta$ use the induction
hypothesis.  Otherwise, $\zeta < \varepsilon$ by the definition of
``$h$ witnesses $\bar\eta \rest u$ is free" and the choice of
$B^\varepsilon_\alpha$ in $\boxplus_{5.1}(c)$.  Now use the choice of
$\bold c^\varepsilon_\alpha(a)$ in Case 2 above.]

So indeed $\boxplus_5$ holds.]
\mn
\begin{enumerate}
\item[$\boxplus_6$]  $\chr(G_M) > \kappa$ if $M \in K^*_{\cA}$.
\end{enumerate}
\mn
Why?  Toward contradiction assume $\bold c:G_M \rightarrow \kappa$ is
a colouring.  For each $\eta \in \cA$ and $\varepsilon < \kappa$ let
$\Lambda_{\eta,\varepsilon} = \{\nu:\nu \in \cA,\nu <_{\cA}
\eta,\nu(\varepsilon) = \eta(\varepsilon)$ and for some $a \in
P^M_\nu$ we have $\bold c(a) = \varepsilon\}$.

Let $\cB_\varepsilon = \{\eta \in \cA:|\Lambda_{\eta,\varepsilon}| <
\kappa\}$.  Now if $\cA \ne \cup\{\cB_\varepsilon:\varepsilon 
< \kappa\}$ then pick any $\eta \in
\cA \backslash \cup\{\cB_\varepsilon:\varepsilon < \kappa\}$ and by
induction on $\varepsilon < \kappa$ choose $\nu_\varepsilon \in
\Lambda_{\eta,\varepsilon} \backslash \{\nu_\zeta:\zeta <
\varepsilon\}$, possible as $\eta \notin \cB_\varepsilon$ 
by the definition of $\cB_\varepsilon$.  By the definition of
$\Lambda_{\eta,\varepsilon}$ there is $a_\varepsilon \in
P^M_{\nu_\varepsilon}$ such that $\bold c(\nu_\varepsilon) =
\varepsilon$.  So as $M \in K^*_{\cA}$ there is $a \in P^M_\eta$ 
such that $\varepsilon < \kappa \Rightarrow
F^M_\varepsilon(a) = a_\varepsilon$, but $\{a,a_\varepsilon\} \in
\edge(G_M)$ hence $\bold c(a) \ne \bold c(a_\varepsilon) =
\varepsilon$ for every $\varepsilon < \kappa$, contradiction.  So $\cA
= \cup\{\cB_\varepsilon:\varepsilon < \kappa\}$.

For each $\varepsilon < \kappa$ we choose $\zeta_\eta < \kappa$ for
$\eta \in \cB_\varepsilon$ by induction on $<_{\cA}$ such
that $\zeta_\eta \notin \{\zeta_\nu:\nu \in
\Lambda_{\eta,\varepsilon} \cap \cB_\varepsilon\}$.  Let
$\cB_{\varepsilon,\zeta} = \{\eta \in \cB_\varepsilon:\zeta_\eta
=\zeta\}$ for $\varepsilon,\zeta < \kappa$ so $\cA =
\cup\{\cB_{\varepsilon,\zeta}:\varepsilon,\zeta < \kappa\}$
 and clearly $\eta \mapsto \eta(\varepsilon)$ is a 
one-to-one function with domain $\cB_{\varepsilon,\zeta}$, 
contradiction to ``$\bar\eta = \bar\eta \rest u_\lambda$ is not weakly free".
\end{PROOF}

\begin{observation}
\label{c6}
1) If $\cA \subseteq {}^\kappa \mu$ and $\eta \ne \nu \in \cA
   \Rightarrow (\forall^\infty \varepsilon < \kappa)(\eta(\varepsilon)
   \ne \nu(\varepsilon))$ then $\cA$ is free iff $\cA$ is weakly
   free.

\noindent
2) The assumptions of \ref{c3}(2) hold \when \,: $\mu \ge \lambda > \kappa$ are
regular, $S \subseteq S^\mu_\kappa$ stationary, $\bar\eta =
   \langle \eta_\delta:\delta \in S\rangle,\eta_\delta$ an increasing
   sequence of ordinals of length $\kappa$ with limit $\delta$ such
   that $u \subseteq [\lambda]^{< \lambda} \Rightarrow \langle
   \Rang(\eta_\delta):\eta \in u\rangle$ has a one-to-one choice function.
\end{observation}

\begin{conclusion}
\label{c12}
Assume that for every graph $G$, if $H \subseteq G \wedge |H| <\lambda
\Rightarrow \chr(H) \le \kappa$ then $\chr(G) \le \kappa$.

\Then \,:
\mn
\begin{enumerate}
\item[$(A)$]  if $\mu > \kappa = \cf(\mu)$ and $\mu \ge \lambda$ then
  $\pp(\mu) = \mu^+$
\sn
\item[$(B)$]  if $\mu > \cf(\mu) \ge \kappa$ and $\mu \ge \lambda$ then
$\pp(\mu) = \mu^+$, i.e. the strong hypothesis
\sn
\item[$(C)$]  if $\kappa = \aleph_0$ then above $\lambda$ the SCH
  holds.
\end{enumerate}
\end{conclusion}

\begin{PROOF}{\ref{c12}}

\noindent
\underline{Clause $(A)$}:  By \ref{c3} and \cite[Ch.II]{Sh:g},
\cite[Ch.IX,\S1]{Sh:g}.
\medskip

\noindent
\underline{Clause $(B)$}:  Follows from (A) by
\cite[Ch.VIII,\S1]{Sh:g}.
\medskip

\noindent
\underline{Clause $(C)$}:  Follows from (B) by \cite[Ch.IX,\S1]{Sh:g}.
\end{PROOF}


\end{document}